\documentclass[12pt,twoside]{article}
\usepackage{latexsym}
\usepackage{amssymb, amsmath}
\usepackage{bm}
\usepackage{fleqn}
\newcommand\chapter{\pagestyle{myheadings}}
\def\cP{{\cal P}}
\def\cI{{\cal I}}
\def\cT{{\cal T}}
\def\cK{{\cal K}}
\def\cO{{\cal O}} 
\def\cN{{\cal N}} 
\def\R{\mathbb{R}}
\def\N{\mathbb{N}}
\def\F{{\sf F}} 
\def\veca{\bm{a}}

\def\vece{\bm{e}}
\def\vech{\bm{h}}
\def\vecn{\bm{n}}
\def\vect{\bm{t}}
\def\vecu{\bm{u}}
\def\vecw{\bm{w}}
\def\vecx{\bm{x}}
\def\vecy{\bm{y}}
\def\vec0{\bm{0}}
\def\vecxi{\bm{\xi}}
\def\veczeta{\bm{\zeta}}
\def\eps{\varepsilon}
\def\ov{\overline} 
\def\DIV{\mathop\mathrm{div}}
\def\disp{\displaystyle } 
\def\<{\langle}
\def\>{\rangle}

\makeatletter 
\long\def\@makefntext#1{\parindent 1em\noindent 
\@hangfrom{\hbox to 1.8em{\hss$^{\@thefnmark}$}}#1}
\makeatother
\def\prf{\noindent{\it Proof}.\ } 
\def\qed{{\hfill\QED}}            
\newcommand{\QED}{\hbox{\rule[0pt]{3pt}{6pt}}}
\def\barint_#1{\mathchoice
  {\mathop{\vrule width 6pt height 3 pt depth -2.5pt
    \kern -8.8pt \intop}\nolimits_{#1}}%
  {\mathop{\vrule width 5pt height 3 pt depth -2.6pt 
    \kern -6.5pt \intop}\nolimits_{#1}}%
  {\mathop{\vrule width 5pt height 3 pt depth -2.6pt
    \kern -6pt \intop}\nolimits_{#1}}%
  {\mathop{\vrule width 5pt height 3 pt depth -2.6pt
    \kern -6pt \intop}\nolimits_{#1}}}
\makeatletter
 
 \@addtoreset{equation}{section}
\makeatother

\newtheorem{Th}{Theorem}[section] 
\newtheorem{Lem}[Th]{Lemma} 
\newtheorem{Prop}[Th]{Proposition}

\newtheorem{Prob}[Th]{Problem}

\begin{document}
\begin{center}
{\large \bf ANALYTICAL AND NUMERICAL ASPECTS ON MOTION OF POLYGONAL CURVES WITH \\ CONSTANT AREA SPEED}
\vskip 5mm
{\large \it Michal Bene\v s, Masato Kimura, and Shigetoshi Yazaki}
\end{center}
\thispagestyle{empty}
\section{\normalsize Introduction}

The first purpose of this paper is to propose a formulation of 
general area-preserving motion of polygonal curves by using a system of ODEs. 
Solution polygonal curves belong to a prescribed polygonal class, 
which is similar to admissible class used in the so-called crystalline curvature flow. 
Actually, if the initial curve is a convex polygon, 
then our polygonal flow is nothing but the crystalline curvature flow. 
However, 
if the initial polygon is not convex and does not belong to any admissible class, 
then the polygonal flow cannot be regarded as a crystalline curvature flow. 
Because the prescribed polygonal class is determined by the initial polygon and 
one can take any polygon as the initial data. 
On the other hand, 
in the framework of the crystalline curvature flow, 
the initial polygon should be taken from the admissible class. 

The second purpose is to discretize 
the ODEs implicitly in time keeping a given constant area speed, 
while the solution polygonal curve exists in the prescribed polygonal class. 

The organization of this paper is as follows. 
In the next section, we will introduce notion of polygonal motion and a polygonal class will be given. 
In section 3, our problem will be formulated and some examples will be given. 
In the last section, we will propose a scheme which achieves the second purpose and 
show convergence of the scheme. 

\section{\normalsize Polygons and polygonal motions}\label{basic} 
\subsection{\normalsize Polygons}\label{polygons} 

We define a set of polygons in $\R^2$: 
\[
\cP:=\{\Gamma;\ \mbox{$\Gamma$ is a polygonal Jordan curve in $\R^2$}\}. 
\]
For $\Gamma\in\cP$, the bounded interior polygonal domain
surrounded by $\Gamma$ is denoted by $\Omega$.
For simplicity, we consider the case that
$\Omega$ is simply connected,
but many of the following arguments are valid 
in other geometrical situations, even in three dimensional case, 
with some minor changes.

Let $\Gamma\in \cP$ be an $N$-polygon.
The $N$ vertices of $\Gamma$ are denoted by
$\vecw_j\in \R^2$
for $j=1, 2, \ldots, N$ counterclockwise,
where $\vecw_0=\vecw_N$ and $\vecw_{N+1}=\vecw_1$.
Hereafter we use the periodic boundary condition 
$\F_0=\F_N$ and $\F_{N+1}=\F_1$ 
for any quantities defined on $N$-polygon. 

The $j$th edge between $\vecw_{j-1}$ and $\vecw_j$ is 
\[
\Gamma_j=\{(1-\theta)\vecw_{j-1}+\theta\vecw_j;\ 0<\theta<1\}
\quad (j=1, 2, \ldots, N),
\]
and their lengths are 
$|\Gamma_j|:=| \vecw_j-\vecw_{j-1}|$.
We define $\chi_j\in L^\infty(\Gamma)$ as
\[
\chi_j(\vecx):=
\left\{
\begin{array}{@{}ll@{}}
1, & \vecx\in\Gamma_j \\
0, & \vecx\in\Gamma\setminus\Gamma_j
\end{array}
\right.
\quad (j=1, 2, \ldots, N),
\]
which is the characteristic function of $\Gamma_j$.

The counterclockwise tangential unit vector
and the outward unit normal are 
denoted by $\vect_j$ and $\vecn_j$, where
$\vect_j=(\vecw_{j}-\vecw_{j-1})/|\Gamma_j|$
and 
$\vecn_j$ is defined such as $\det(\vecn_j, \vect_j)=1$. 
The outer angle at the vertex $\vecw_j$ is denoted
by $\varphi_j\in (-\pi,\pi)\setminus \{0\}$. 

We remark that
\[
\cos \varphi_j= \vect_{j+1}\cdot\vect_j=\vecn_{j+1}\cdot\vecn_j .
\]
We define the height of $\Gamma_j$ from the origin
$h_j:=\vecw_j\cdot\vecn_j$. 
They satisfy the equalities
\begin{equation}\label{aba}
|\Gamma_j|=a_{j-1}h_{j-1}
+b_jh_j+a_jh_{j+1}
\quad (j=1, 2, \ldots, N),
\end{equation}
where 
$a_j:=(\sin \varphi_j)^{-1}$ and $b_j:=-\cot \varphi_{j-1}-\cot \varphi_j$. 
This equality can be checked from the fact that
the straight line 
including $\Gamma_j$ is expressed by the equation
$\vecn_j\cdot\vecx=h_j$.
The total length of $\Gamma$ is given by
\begin{equation}\label{Glength}
|\Gamma|:=\sum_{j=1}^N |\Gamma_j|
=\sum_{j=1}^N (a_j+b_j+a_{j-1})h_j
=\sum_{j=1}^N \eta_j h_j,
\end{equation}
where
$\eta_j:=a_j+b_j+a_{j-1}=\tan(\varphi_j/2)+\tan(\varphi_{j-1}/2)$. 

The area of interior domain $\Omega$ is denoted by $|\Omega|$,
which is given by 
\begin{equation}\label{area0}
|\Omega|=\frac{1}{2}\sum_{j=1}^N|\Gamma_j|h_j.
\end{equation}

The above symbols are also written as
$\Omega=\Omega(\Gamma)$, 
$\vecn_j=\vecn_j(\Gamma)$ and $h_j=h_j(\Gamma)$ etc., 
if we need to distinguish 
from quantities of the other polygons.

\subsection{\normalsize Motion of polygons and Lipschitz mappings}\label{motion} 
\setcounter{equation}{0}

We consider a moving polygon $\Gamma(t)\in\cP$,
where the parameter $t$ (we call $t$ time)
belongs to an interval $\cI\subset\R$. 
For $k\in\N\cup\{0\}$, 
we call a moving polygon $\Gamma(t)$ belongs to
$C^k$-class on $\cI$, if the number of edges of $\Gamma(t)$
does not change in time and $\vecw_j\in C^k(\cI, \R^2)$ for all $j=1, 2, \ldots, N$.

If $k\geq 1$, we can define the normal velocity at 
$\vecx\in\Gamma_j(t)$ which is the $j$th edge of $\Gamma(t)$. 
We suppose $\vecx^*\in\Gamma_j(t^*)$ and
$\vecx^*=(1-\theta)\vecw_{j-1}(t^*)+\theta\vecw_j(t^*)$ for some $\theta\in(0,1)$, 
and define $\vecx(\theta,t):=(1-\theta)\vecw_{j-1}(t)+\theta\vecw_j(t)\in\Gamma(t)$. 
Then the outward normal velocity of $\Gamma_j(t^*)$ at $\vecx^*$ is defined by
\[
V_j(\vecx^*, t^*):=\dot{\vecx}(\theta, t^*)\cdot\vecn_j(t^*)
=(1-\theta)\dot{\vecw}_{j-1}(t^*)\cdot\vecn_j(t^*)
+\theta\dot{\vecw}_j(t^*)\cdot\vecn_j(t^*). 
\]
Here and hereafter, 
the (partial) derivative of $\F$ with respect to $t$ is denoted by $\dot{\F}$. 
We remark that $V_j(\cdot,t)$ is a linear function on each $\Gamma_j(t)$. 
We define the normal velocity of $\Gamma(t)$ by
\[
V(\cdot, t):=\sum_{j=1}^N V_j(\cdot, t)\chi_j(\cdot, t)\in L^\infty(\Gamma(t)),
\]
where $\chi_j(\cdot, t)\in L^\infty(\Gamma(t))$ is
the characteristic function of $\Gamma_j(t)$. 

For a $C^k$-class moving $N$-polygon $\Gamma(t)$ 
$(t\in \cI)$ with its interior domain $\Omega(t)$, 
we construct Lipschitz mappings smoothly parametrized by $t$ from a fixed polygonal domain. 

We fix $t^*\in\cI$ and define $\Gamma^*:=\Gamma (t^*)$ and
$\Omega^*:=\Omega (t^*)$.
We also choose another domain $Q$ with 
$\ov{\cup_{t\in\cI}\Omega (t)}\subset Q$ and fix it. 
Our aim of this subsection is to construct Lipschitz mappings 
$\Phi(t)=\Phi(\cdot, t)$ from $\Omega^*$ to $\Omega(t)$ smoothly parametrized by $t$. 

\begin{Prop}\label{propmapping}
Under the above condition, 
there exists $\eps>0$, $\cI^*:=\cI\cap(t^*-\eps,\ t^*+\eps)$, 
and $\Phi\in C^k(\cI^*, W^{1, \infty}(Q, \R^2))$, 
and they satisfy the following conditions. 
\begin{enumerate}
 \item $\Phi(t)$ is a bi-Lipschitz transform from $\ov{Q}$ onto itself, 
	i.e. $\Phi(t)$ is bijective from $\ov{Q}$ onto itself 
	and $\Phi(t)$ and $\Phi(t)^{-1}$ are both Lipschitz continuous on $\ov{Q}$. 
 \item $\Phi(\vecx, t)=\vecx$ for $\vecx$ in a neighborhood of $\partial Q$ for $t\in\cI^*$. 
 \item $\Phi(\ov{\Omega^*}, t)=\ov{\Omega(t)}$ for $t\in\cI^*$, 
	and $\Phi(t)$ is an affine map from each edge $\Gamma_j^*:=\Gamma_j(t^*)$ onto $\Gamma_j(t)$ 
	with $\Phi(\vecw_j(t^*),t)=\vecw_j(t)$. 
\end{enumerate}
\end{Prop}
\prf 
Without loss of generality,
we assume that $Q$ is a bounded polygonal domain, too.
We consider a triangulation $\cT$ of $\Omega^*$ and $Q$. 
Namely $\cT$ is a collection of triangular subatomics of $Q$ with 
\[
\mbox{
$\disp\ov{Q}=\bigcup_{K\in\cT}\ov{K}$, 
\ and\ 
$\disp K\cap K'=\emptyset$ 
\ if\ 
$\disp K, K'\in\cT$, 
$\disp K\neq K'$}, 
\]
and $\ov{K}\cap\ov{K'}$ is either the empty set, 
a common vertex or a common edge of $K$ and $K'\in\cT$, 
and there is a subset $\cT_0\subset \cT$ such that 
\[
\mbox{
$\disp\ov{\Omega^*}=\bigcup_{K\in\cT_0}\ov{K}$, 
\ and\ 
$\disp\Gamma^*\cap \cN =\{\vecw_j^*\}_{j=1}^N$}, 
\]
where $\cN$ denotes the set of all vertices of triangles in $\cT$ 
and $\vecw_j^*:=\vecw_j(t^*)$. 
We also suppose that there does not exist any $K\in\cT$ 
with $\ov{K}\cap\Gamma^*\neq\emptyset$ and $\ov{K}\cap\partial Q \neq\emptyset$. 

We assume that $\Phi(\vecx, t)$ has the following form: 
\begin{equation}\label{PhiM}
\Phi(\vecx,t)=A_K(t)
\left(
\begin{array}{@{}c@{}}
\vecx\\ 1
\end{array}
\right)
\quad (\vecx\in K\in\cT),
\end{equation}
where 
$A_K(t)$ is a $2\times 3$ matrix depending on $K\in\cT$ and $t$. 
To determine $A_K(t)$, we suppose the condition:
\begin{equation}\label{Phixit}
\Phi(\vecx, t)=
\left\{
\begin{array}{@{}ll}
\vecx & \mbox{if $\vecx\in\cN\setminus\Gamma^*$}, \\
\vecw_j(t) & \mbox{if $\vecx=\vecw_j^*\in\cN\cap\Gamma^*$}.
\end{array}
\right.
\end{equation}
For sufficiently small $\eps>0$, for $t\in\cI^*$, $\Phi(\vecx, t)$ is uniquely determined 
by the conditions (\ref{PhiM}) and (\ref{Phixit}). 
It is also clear that $\Phi(t)=\Phi(\cdot, t)\in W^{1,\infty}(Q,\R^2)$ 
is bijective from $Q$ onto itself and $\Phi(\vecx, t)=\vecx$ for $\vecx\in K$ 
if $\ov{K}\cap\Gamma^*=\emptyset$. 

Let us fix $K\in\cT$ (with $\ov{K}\cap\Gamma^*\neq\emptyset$) 
and let $\vecx_1$, $\vecx_2$, $\vecx_3$ be the vertices of $K$. 
Then, from (\ref{Phixit}), $\vecy_l(t):=\Phi(\vecx_l, t)$ 
satisfies the condition $\vecy_l\in C^k(\cI^*,\R^2)$ for $l=1, 2, 3$. 
Since 
\[
\vecy_l=A_K(t)
\left(
\begin{array}{@{}c@{}}
\vecx_l\\ 1
\end{array}
\right)
\quad (l=1, 2, 3), 
\]
we have
\[
A_K(t)
=(\vecy_1(t), \vecy_2(t), \vecy_3(t))
\left(
\begin{array}{@{}ccc@{}}
\vecx_1 & \vecx_2 & \vecx_3  \\
1 & 1 & 1
\end{array}
\right)^{-1}. 
\]
From this expression, we obtain $A_K\in C^k(\cI^*, \R^{2\times 3})$ 
and hence $\Phi\in C^k(\cI^*, W^{1,\infty}(Q,\R^2))$ follows. 
\qed

As an application of this proposition, 
we have the following theorem which is well-known 
in the case $\Gamma(t)$ is a smooth Jordan curve. 
\begin{Th}\label{ddtint}
Let $\Gamma(t)\in \cP$ be a $C^1$-class moving $N$-polygon 
on an interval $t\in \cI$ with its interior domain $\Omega(t)$. 
For all $\phi\in C^1(\R^2\times\cI)$, 
the map 
$[t\mapsto\int_{\Omega(t)}\phi(\vecx, t)\,d\vecx]$ 
belongs to $C^1(\cI)$ and 
\[
\frac{d}{dt}\int_{\Omega(t)}\phi(\vecx,t)\,d\vecx 
=\int_{\Omega(t)}\dot{\phi}(\vecx,t)\,d\vecx 
+\int_{\Gamma(t)}\phi(\vecx,t)V(\vecx,t)\,ds
\]
holds. 
\end{Th}
\prf 
We define $f(t):=\int_{\Omega(t)}\phi(\vecx,t)\,d\vecx$.
Under the setting of Proposition~\ref{propmapping}, 
we show $f\in C^1(\cI^*)$ and calculate $\dot{f}(t^*)$. 
For $t\in \cI^*$, we have
\[
f(t)=\int_{\Omega^*}
\phi(\Phi(\vecx, t), t)J(\vecx, t)\,d\vecx,
\]
where $J(\vecx,t)$ is the Jacobian defined by 
$J(\vecx,t):=\det(\nabla_{\vecx}\Phi^{\rm T}(\vecx, t))$. 
We remark that $J\in C^1(\cI^*, L^\infty(Q))$ 
and $\dot{J}(\vecx, t^*)={\DIV}_{\vecx}\dot{\Phi}(\vecx, t^*)$ 
since $\Phi\in C^1(\cI^*, W^{1, \infty}(Q, \R^2))$ and $\Phi(\vecx, t^*)=\vecx$. 
Hence, we obtain 
\begin{eqnarray*}
\dot{f}(t^*)&=&
\int_{\Omega^*}
(\nabla_{\vecx}\phi(\vecx, t^*)\cdot\dot{\Phi}(\vecx, t^*)+\dot{\phi}(\vecx,t^*))\,d\vecx
+\int_{\Omega^*}\phi(\vecx,t^*)\dot{J}(\vecx, t^*)\,d\vecx \\
&=&
\int_{\Omega^*}
{\DIV}_{\vecx}(\phi(\vecx, t^*)\dot{\Phi}(\vecx, t^*))\,d\vecx
+\int_{\Omega^*}\dot{\phi}(\vecx, t^*)\,d\vecx \\
&=&
\int_{\Gamma^*}\phi(\vecx, t^*)V(\vecx, t^*)\,ds
+\int_{\Omega^*}\dot{\phi}(\vecx, t^*)\,d\vecx
\end{eqnarray*}
\qed

In the case where $\phi\equiv 1$, 
we obtain the following formula for $C^1$-class moving polygon particularly. 
\begin{equation}\label{area}
\frac{d}{dt}|\Omega(t)|=\int_{\Gamma(t)}V(\vecx, t)\,ds,
\end{equation}
where $|\Omega(t)|$ stands for the area of $\Omega(t)$. 

We remark that Proposition \ref{propmapping} and Theorem \ref{ddtint} 
with their proofs are valid even in the three dimensional case.

\subsection{\normalsize Polygonal motion}\label{polygonal} 

For two polygons $\Gamma$ and $\Sigma\in\cP$, 
we define an equivalence relation $\Gamma\sim\Sigma$. 
We say $\Gamma\sim\Sigma$, 
if their numbers of edges are same (let it be $N$) and 
$\vecn_j(\Gamma)=\vecn_j(\Sigma)$ for all $j=1, 2, \ldots, N$ 
after choosing suitable counterclockwise numbering for $\Gamma$ and $\Sigma$. 
The equivalence class of $\Gamma\in\cP$ is denoted by 
$\cP[\Gamma]:=\{\Sigma\in\cP;\ \Sigma\sim\Gamma\}$. 

We fix an $N$-polygon $\Gamma^*\in\cP$ and let $\cP^*:=\cP[\Gamma^*]$. 
For $\Gamma$ and $\Sigma$ in $\cP^*$, 
we define the distance between them by 
\[
d(\Gamma, \Sigma):=\max_{j=1, 2, \ldots, N}|h_j(\Gamma)-h_j(\Sigma)|. 
\]
Then, it is clear that $(\cP^*, d)$ becomes a metric space since it is isometrically embedded in 
$\R^N$ equipped with maximum norm $|\cdot|_\infty$ by the height function $\vech$: 
$\cP^*\ni\Gamma\mapsto 
\vech(\Gamma)=(h_1(\Gamma), h_2(\Gamma), \ldots, h_N(\Gamma))^{\rm T}\in\R^N$. 
The following proposition is clear. 
\begin{Prop}
The set $\vech(\cP^*)$ is an open subset of $\R^N$.
\end{Prop}

For any $\Gamma^0$ and $\Gamma^1\in\cP^*$ and for $\theta\in [0,1]$, we define 
\[
\vech^\theta:=(1-\theta)\vech (\Gamma^0)+\theta \vech (\Gamma^1)\in \R^N. 
\]
If there exists $\Gamma^\theta\in\cP^*$ with $\vech(\Gamma^\theta)=\vech^\theta$,
$\Gamma^\theta$ is called $\theta$-interpolation of $\Gamma^0$ and $\Gamma^1$.
The $\theta$-interpolation of $\Gamma^0\in\cP^*$
and $\Gamma^1\in\cP^*$ is denoted by 
$(1-\theta)\Gamma^0+\theta\Gamma^1:=\Gamma^\theta\in\cP^*$.
We remark that it satisfies 
\[
|\Gamma^\theta_j|\geq \min\{|\Gamma^0_j|,\ |\Gamma^1_j|\}.
\]

For $\Gamma\in\cP^*$ and $\eps>0$,
$\eps$-ball in $\cP^*=\cP[\Gamma]$
with center $\Gamma$ is denoted by
\[
B(\Gamma, \eps):=
\{\Sigma\in\cP[\Gamma];\ d(\Sigma, \Gamma)<\eps\}.
\]
For an open set $\cO\subset\cP^*$ and $\Gamma\in\cO$, 
we define a positive number $\rho(\Gamma, \cO)>0$ as 
\[
\rho(\Gamma, \cO):=
\inf \{|\veca-\vech(\Gamma)|_\infty;\ 
\veca\in\R^N\setminus \vech(\cO)\}.
\]
We remark that 
$\rho(\cdot, \cO)$ is Lipschitz continuous with Lipschitz constant $1$: 
\[
|\rho(\Gamma,\cO)-\rho(\Sigma,\cO)|
\leq
d(\Gamma, \Sigma)
\quad (\Gamma, \Sigma\in\cO).
\]
For a compact set $\cK\subset\cO$, we also define
\[
\rho(\cK, \cO):=\min_{\Gamma\in\cK}\rho(\Gamma,\cO).
\]

Let $a_j^*:=a_j[\Gamma^*]$ and $b_j^*:=b_j[\Gamma^*]$. 
Then, from the formula (\ref{aba}), we obtain
$$
\begin{array}{@{}l}\disp
||\Gamma_j|-|\Sigma_j|| \\[5pt]\disp \qquad
=|a_{j-1}^*(h_{j-1}[\Gamma]-h_{j-1}[\Sigma])
+b_{j}^*(h_{j}[\Gamma]-h_{j}[\Sigma])
+a_{j}^*(h_{j+1}[\Gamma]-h_{j+1}[\Sigma])| \\[5pt]\disp\qquad
\leq C^* d(\Gamma, \Sigma)
\quad (j=1, 2, \ldots, N),
\end{array}
$$
where we define
\begin{equation}\label{C*}
C^*:=\max_{l=1, 2, \ldots, N}\{|a_{l-1}^*|+|b_l^*|+|a_l^*|\}.
\end{equation}
For a compact set $\cK\subset \cP^*$, we define
\[
\sigma(\cK):=\min\{|\Gamma_j|;\ \Gamma\in \cK,\ j=1, 2, \ldots, N\}>0. 
\]

We consider a $C^k$-class moving polygon 
$\Gamma (t)\in\cP^*$ $(t\in\cI)$. 
We call it {\bf polygonal motion in $\bm{\cP^*}$} in this paper. 
We remark that a polygonal motion 
$\Gamma(t)\in\cP^*$ $(t\in\cI)$ belongs to $C^k$-class if and only if 
$h_j\in C^k(\cI)$ for $j=1, 2, \ldots, N$. 
If $\Gamma(t)$ is a $C^1$-class polygonal motion, 
its normal velocity $V_j$ of $\Gamma_j(t)$ 
is a constant on each $\Gamma_j(t)$ and it is given by $V_j(t)=\dot{h}_j(t)$. 
The formula (\ref{area}) is written in the form: 
\begin{equation}\label{area2}
\frac{d}{dt}|\Omega(t)|
=\sum_{j=1}^N|\Gamma_j(t)|V_j(t). 
\end{equation}
 
We fix an equivalence class $\cP^*$ of polygons 
and let its $j$th outward unit normal be $\vecn_j$ and outer angle $\varphi_j$. 
For $\Gamma\in\cP^*$, the {\bf polygonal curvature} 
$\kappa_j$ of $\Gamma_j$ is defined by
\[
\kappa_j:=\frac{\eta_j}{|\Gamma_j|}. 
\]
We also define the polygonal curvature of $\Gamma$ by 
\[
\kappa:=\sum_{j=1}^N
\kappa_j\chi_j\in L^\infty(\Gamma). 
\]
The reason why this is called ``curvature'' is shown by the following proposition. 
\begin{Prop}
Let $\Gamma(t)$ $(t\in\cI)$ be a $C^1$-class polygonal motion in $\cP^*$. Then 
\[
\frac{d}{dt}|\Gamma(t)|
=\sum_{j=1}^N|\Gamma_j(t)|\kappa_j(t)V_j(t)
=\int_{\Gamma(t)}\kappa(\vecx,t)V(\vecx,t)\,ds. 
\]
\end{Prop}
\prf 
We obtain 
\[
\frac{d}{dt}|\Gamma(t)|
=\frac{d}{dt}\sum_{j=1}^{N}\eta_jh_j(t)
=\sum_{j=1}^{N}\eta_jV_j(t)
=\sum_{j=1}^N|\Gamma_j(t)|\kappa_j(t)V_j(t), 
\]
from the formula (\ref{Glength}). 
\qed
\section{\normalsize Initial value problem of polygonal motion}\label{IVP} 
\subsection{\normalsize General polygonal motion problems}\label{general} 

We fix an equivalence class of $N$-polygons $\cP^*$ as in \S\ref{polygonal}.
For an open set $\cO\subset\cP^*$ and $T_*\in(0, \infty]$, 
let $F$ be a given continuous function from 
$\cO\times[0,T_*)$ to $\R^N$ with the local Lipschitz property: 
For arbitrary compact set $\cK\subset\cO$ and $T\in (0, T_*)$, 
there exists $L(\cK, T)>0$ such that 
\begin{equation}\label{Lip-cond}
|F(\Gamma,t)-F(\Sigma,t)|_\infty\leq L(\cK, T)\, d(\Gamma, \Sigma)
\quad (\Gamma, \Sigma\in \cK,\ t\in [0,T]). 
\end{equation}
Under the condition (\ref{Lip-cond}), 
for a compact set $\cK\subset\cO$ and $T\in (0, T_*)$, 
we also define
\[
M(\cK,T):=
\max\{|F(\Gamma,t)|_\infty;\ \Gamma\in \cK,\ t\in [0,T]\}>0. 
\]

We consider the following initial value problem of polygonal motion. 
\begin{Prob}\label{PF}
For a given $N$-polygon $\Gamma^*\in\cO$, 
find a $C^1$-class polygonal motion $\Gamma(t)\in\cO$
$(0\leq t\leq T<T_*)$ such that
\[
\left\{
\begin{array}{@{}l}\disp
V_j(t)=F_j(\Gamma(t),t)
\quad (t\in [0,T],\ j=1, 2, \ldots, N)\\[5pt]\disp
\Gamma(0)=\Gamma^*.
\end{array}
\right.
\]
\end{Prob}
Under the Lipschitz condition (\ref{Lip-cond}), 
it is clear that there exists a local solution $\Gamma(t)$ in a short time interval $[0,T]$, 
since Problem \ref{PF} can be expressed by 
an initial value problem of an ordinary differential equations 
for $\vech(t)=(h_1(t), h_2(t), \ldots, h_N(t))^{\rm T}$. 

We often assume the following condition for $F_j$:
\begin{equation}\label{areaspeed}
\sum_{j=1}^N|\Gamma_j|F_j(\Gamma,t)=\mu
\quad (\Gamma\in\cO, \ t\in [0,T_*)), 
\end{equation}
where $\mu$ is a fixed real number. 
Under the assumption (\ref{areaspeed}), 
from the formula (\ref{area2}), 
any solution $\Gamma(t)$ to Problem \ref{PF} has the property of 
the following constant area speed ({\bf CAS} in short): 
\[
\frac{d}{dt}|\Omega(t)|=\mu.
\]

The polygonal flow is regarded as the crystalline curvature flow if 
the initial polygon $\Gamma^*$ is convex 
(as mentioned in introduction, if $\Gamma^*$ is not convex, 
then the polygonal flow is different from the crystalline curvature flow). 
There are many articles about the crystalline curvature flow and 
asymptotic behavior of solutions 
\cite{Andrews2002, GigaG2000, IshiiS1999, IshiwataUYY2004, IshiwataY2003, Yazaki2002a, Yazaki2002b, Yazaki2007a, Yazaki2007b}, etc., 
which started from the pioneer works \cite{AngenentG1989} and \cite{Taylor1993}. 

\subsection{\normalsize Examples of problems of polygonal motion}\label{ecmp} 

In this section, we give some basic examples of polygonal motions which are 
nice polygonal analogues of corresponding smooth moving boundary problems. 

\begin{Prob}[polygonal curvature flow]\label{PCF} 
For a given $N$-polygon $\Gamma^*\in\cP^*$, 
find a $C^1$-class family of $N$-polygons $\bigcup_{0\leq t\leq T}\Gamma(t)\subset\cP^*$ $(T<T_*)$ satisfying
\[
\left\{
\begin{array}{@{}l}\disp
V_j(t)=-\kappa_j(t)
\quad (t\in [0,T],\ j=1, 2, \ldots, N), \\[5pt]\disp
\Gamma(0)=\Gamma^*.
\end{array}
\right.
\]
\end{Prob}
The solution has CAS property with $\mu=-2\sum_{j=1}^N\tan(\varphi_j/2)$: 
\[
\frac{d}{dt}|\Omega(t)|
=-\sum_{j=1}^N\kappa_j(t)|\Gamma_j(t)|
=-\sum_{j=1}^N\eta_j
=-2\sum_{j=1}^N\tan\frac{\varphi_j}{2}
=const.
\]

\begin{Prob}[area-preserving polygonal curvature flow]\label{AP-PCF} 
\quad
For a given $N$-polygon $\Gamma^*\in\cP^*$, 
find a $C^1$-class family of $N$-polygons $\bigcup_{0\leq t\leq T}\Gamma(t)\subset\cP^*$ $(T<T_*)$ satisfying
\[
\left\{
\begin{array}{@{}l}\disp
V_j(t)=\<\kappa(\cdot, t)\>-\kappa_j(t)
\quad (t\in [0,T],\ j=1, 2, \ldots, N), \\[5pt]\disp
\Gamma(0)=\Gamma^*.  
\end{array}
\right.
\]
Here $\<\kappa(\cdot, t)\>$ is the mean value of $\kappa$ on $\Gamma(t)${\rm :} 
\[
\<\kappa(\cdot, t)\>
=\frac{1}{|\Gamma(t)|}\int_{\Gamma(t)}\kappa(\vecx, t)\,ds
=\frac{\sum_{i=1}^N\eta_i}{|\Gamma(t)|}
=\frac{2\sum_{i=1}^N\tan(\varphi_i/2)}{|\Gamma(t)|}. 
\]
\end{Prob}
The solution has CAS property with $\mu=0$:
\[
\frac{d}{dt}|\Omega(t)|
=\<\kappa(\cdot, t)\>|\Gamma(t)|-\int_{\Gamma(t)}\kappa(\vecx, t)\,ds
=0. 
\]

In what follows, 
the mean value of $\F$ on the edge $\Gamma_j$ is denoted by 
\[
\<\F\>_j
:=\frac{1}{|\Gamma_j|}\int_{\Gamma_j}\F(\vecx)\,ds. 
\]
Let $G$ be a bounded Lipschitz domain in $\R^2$. 
We define
\[
\cO_G:=\{\Gamma\in\cP^*;\ 
\Omega(\Gamma)\supset\ov{G}\}. 
\]

\begin{Prob}[area-preserving polygonal advected flow]\label{AP-AF} 
Let us consider a divergence free vector field 
$\vecu\in C^1(\R^2\setminus G;\ \R^2)$ with 
$\DIV\vecu=0$ in $\R^2\setminus G$. 
For a given $N$-polygon $\Gamma^*\in\cO_G$, 
find a $C^1$-class family of $N$-polygons $\bigcup_{0\leq t\leq T}\Gamma(t)\subset\cO_G$ $(T<T_*)$ satisfying 
\[
\left\{
\begin{array}{@{}l}\disp
V_j(t)=\<\vecu\cdot\vecn\>_j
\quad (t\in [0,T],\ j=1, 2, \ldots, N), \\[5pt]\disp
\Gamma(0)=\Gamma^*.
\end{array}
\right.
\]
\end{Prob}
The solution has CAS property with $\disp\mu=\int_{\partial G}\vecn\cdot\vecu\,ds$: 
$$
\frac{d}{dt}|\Omega(t)|
=-\int_{\Gamma(t)}\vecn_j\cdot\vecu\,ds
=-\int_{\Omega(t)\setminus\ov{G}}\DIV\vecu\,d\vecx
+\int_{\partial G}\vecn\cdot\vecu\,ds
=\int_{\partial G}\vecn\cdot\vecu\,ds, 
$$
where 
$\vecn_j$ and $\vecn$ are the outward unit normal vector to $\partial(\Omega(t)\setminus \ov{G})$. 

\begin{Prob}[polygonal Hele-Shaw flow]\label{PHF}
For a given $N$-polygon $\Gamma^*\in\cO_G$ and 
a function $b$ defined on $\partial G\times [0,T]$, 
find a $C^1$-class family of $N$-polygons $\bigcup_{0\leq t\leq T}\Gamma(t)\subset\cO_G$ $(T<T_*)$ satisfying 
\[
\left\{
\begin{array}{@{}l@{}l}\disp
V_j(t)=-\vecn_j\cdot\nabla p(\vecx,t)
&\quad (\vecx\in\Gamma_j(t),\ t\in [0,T],\ j=1, 2, \ldots, N), \\[5pt]\disp
\Delta p(\vecx,t)=0
&\quad (\vecx\in\Omega(t)\setminus \ov{G},\ t\in [0,T]), \\[5pt]\disp
\<p(\cdot, t)\>_j=\kappa_j(t)
&\quad (t\in [0,T],\ j=1, 2, \ldots, N), \\[5pt]\disp
\frac{\partial p}{\partial\vecn}(\vecx,t)=b(\vecx,t)
&\quad (\vecx\in\partial G,\ t\in [0,T]), \\[5pt]\disp
\Gamma(0)=\Gamma^*.
\end{array}
\right.
\]
\end{Prob}
Here $\partial p/\partial\vecn=\nabla_{\vecx}p\cdot\vecn$ and 
$\vecn_j$ and $\vecn$ are the outward unit normal vector to $\partial(\Omega(t)\setminus \ov{G})$. 
The solution has a given area speed property: 
$$
\frac{d}{dt}|\Omega(t)|
=-\int_{\Gamma(t)}\vecn_j\cdot\nabla p\,ds
=-\int_{\Omega(t)\setminus\ov{G}}\Delta p\,d\vecx
+\int_{\partial G}\frac{\partial p}{\partial\vecn}(\vecx,t)\,ds
=\int_{\partial G}b(\vecx,t)\,ds. 
$$
If $b(\vecx,t)\equiv b_0$ for a given constant $b_0$, 
then the solution has CAS property with $\mu=|\partial G|b_0$. 

\section{\normalsize Numerical schemes}\label{NS} 
\subsection{\normalsize Notation}\label{NNS} 

In \S\ref{NS}, we consider time discretization of Problem \ref{PF} with the following notation. 
The discrete time steps are denoted by 
$0=t_0<t_1<t_2<\cdots<t_{\bar{m}}\leq T$. 
The step size which may be nonuniform and their maximum size are defined by 
\[
\tau_m:=t_{m+1}-t_m
\quad (m=0,1,\cdots,\bar{m}-1), 
\quad \tau:=\max_{0\leq m <\bar{m}}\tau_m. 
\]
Approximate solution of $\Gamma(t_m)$ is denoted by $\Gamma^m\in\cP^*$. 
Quantities of the polygon $\Gamma^m$ are denoted by $h_j^m:=h_j (\Gamma^m)$, 
and $\kappa_j^m:=\kappa_j(\Gamma^m)$, etc. 
We define $e_j^m:=h_j(t_m)-h_j^m$ and $\vece^m:=(e_1^m, e_2^m, \ldots, e_N^m)^{\rm T}\in\R^N$. 
Then we have $d(\Gamma(t_m), \Gamma^m)=|\vece^m|_\infty$. 

The discrete normal velocity $V_j^m$, 
which is an approximation of $V_j(t_m)=\dot{h}_j(t_m)$, is defined by 
\[
V_j^m:=\frac{h_j^{m+1}-h_j^m}{\tau_m}
\quad (m=0, 1, \ldots, \bar{m}-1). 
\]
Corresponding to the formula (\ref{area2}), the following formula holds. 
\begin{equation}\label{daspeed}
\frac{|\Omega^{m+1}|-|\Omega^m|}{\tau_m}
=\sum_{j=1}^N\frac{|\Gamma_j^m| +|\Gamma_j^{m+1}|}{2}V_j^m. 
\end{equation}
This has a form of sum of areas of $N$ trapezoids and is actually derived from (\ref{area0}) as follows: 
\begin{eqnarray*}
|\Omega^{m+1}|-|\Omega^m|
&=&
\frac{1}{2}\sum_{j=1}^N
\left( 
|\Gamma_j^{m+1}|h_j^{m+1}-|\Gamma_j^m|h_j^m
\right)\\
&=&
\frac{1}{2}\sum_{j=1}^N
\left\{
(|\Gamma_j^{m+1}|+|\Gamma_j^m|)
(h_j^{m+1}-h_j^m)
+|\Gamma_j^{m+1}|h_j^m
-|\Gamma_j^m|h_j^{m+1}
\right\}\\
&=&
\frac{\tau_m}{2}\sum_{j=1}^N
(|\Gamma_j^{m+1}|+|\Gamma_j^m|)
V_j^m
+
\frac{1}{2}\sum_{j=1}^N
\left(
|\Gamma_j^{m+1}|h_j^m
-|\Gamma_j^m|h_j^{m+1}
\right) ,
\end{eqnarray*}
where the last sum is equal to zero due to the equality (\ref{aba}).

In the following subsections, 
we suppose that there exists a unique solution $\Gamma(t)$ for $0\leq t\leq T<T_*$ 
to Problem \ref{PF} under the condition (\ref{Lip-cond}), 
and that discrete time steps $0=t_0<t_1<t_2<\cdots<t_{\bar{m}}\leq T$ 
are given a priori such as the uniform time stepping $t_m=m\tau$. 
It is, however, possible to apply any a posteriori adaptive time step control scheme. 

%

\subsection{\normalsize Euler scheme}\label{ES} 

We consider the following Euler scheme to discretized Problem \ref{PF}.
\begin{Prob}\label{es}
For a given $N$-polygon $\Gamma_*\in\cO$ 
and time steps $0=t_0<t_1<t_2<\cdots<t_{\bar{m}}\leq T$, 
find polygons $\Gamma^m\in\cO$ $(m=1, 2, \ldots, \bar{m})$ such that 
\[
\left\{
\begin{array}{@{}l}\disp
V_j^m=F_j(\Gamma^m, t_m)
\quad
(m=0, 1, 2, \ldots, \bar{m}-1,\ j=1, 2, \ldots, N), \\[5pt]\disp
\Gamma^0=\Gamma_*. 
\end{array}
\right.
\]
\end{Prob}
\begin{Th}
We suppose the condition (\ref{Lip-cond}) 
and that $\{\Gamma(t)\}_{0\leq t\leq T}$ 
be a $C^{k+1}$-class solution of Problem {\rm \ref{PF}} 
for $k=0$ or $1$. 
There exists $\delta^*>0$, $\tau^*>0$, $C>0$ and 
a non-decreasing function $\omega (a)>0$ with 
\begin{equation}\label{decayomega-es}
\omega(a)=\left\{
\begin{array}{@{}ll}\disp
o(1) & \mbox{if $k=0$,}\\[5pt]\disp
O(a) & \mbox{if $k=1$,}
\end{array}
\right.
\quad \mbox{as}\ a\downarrow 0,
\end{equation}
such that, if $d(\Gamma^*, \Gamma^0)\leq \delta^*$ and $\tau\leq \tau^*$, 
then $\Gamma^m\in\cO$ $(m=1, 2, \ldots, \bar{m})$ is determined by
the Euler scheme (Problem~\ref{es}) and  
satisfies the estimate
\[
\max_{0\leq m\leq \bar{m}}
d(\Gamma(t_m), \Gamma^m)\leq 
\omega(\tau)+C d(\Gamma(0), \Gamma^0).
\]
\end{Th}
The proof is similar to the one of Theorem \ref{Th:convergence_implicit_2nd_order_scheme}. 

\subsection{\normalsize Second order implicit scheme}\label{2IS} 

We consider the following implicit scheme for Problem \ref{PF}.
\begin{Prob}\label{2is}
For a given $N$-polygon $\Gamma_*\in\cO$ 
and time steps $0=t_0<t_1<t_2<\cdots<t_{\bar{m}}\leq T$, 
find polygons $\Gamma^m\in\cO$
$(m=1, 2, \ldots, \bar{m})$ such that
\[
\left\{
\begin{array}{@{}l}\disp
V_j^m=F_j(\Gamma^{m+1/2}, t_{m+1/2})
\quad (m=0, 1, 2, \ldots, \bar{m}-1,\ j=1, 2, \ldots, N), \\[5pt]
\Gamma^0=\Gamma_*, 
\end{array}
\right.
\]
where $\Gamma^{m+1/2}$ and $t_{m+1/2}$ are the $1/2$-interpolations{\rm :}
\[
\Gamma^{m+1/2}:=\frac{\Gamma^m+\Gamma^{m+1}}{2}\in\cP^*, 
\quad 
t_{m+1/2}:=\frac{t_m+t_{m+1}}{2}=t_m+\frac{\tau_m}{2}.
\]
\end{Prob}
This is a generalized version of \cite{UY2004} for area-preserving crystalline curvature flow. 
\begin{Prop}
We suppose the constant speed condition {\rm (\ref{areaspeed})}. 
Let $\Gamma^m\in\cO$ $(m=1, 2, \ldots, \bar{m})$ be a solution of Problem {\rm \ref{2is}}. 
Then it satisfies 
\[
|\Omega^{m+1}|=|\Omega^m|+\mu \tau_m
\quad (m=0, 1, \ldots, \bar{m}-1). 
\]
In other words, $|\Omega^m|=|\Omega(t_m)|$ holds 
if the exact solution $\Omega(t)$ of Problem {\rm \ref{PF}} exists.
\end{Prop}
\prf 
Since 
$|\Gamma^{m+1/2}_j|=(|\Gamma^m_j|+|\Gamma^{m+1}_j|)/2$, 
we have
\[
\frac{|\Omega^{m+1}|-|\Omega^m|}{\tau_m}
=\sum_{j=1}^N|\Gamma^{m+1/2}_j|F_j(\Gamma^{m+1/2}, t_{m+1/2})
=\mu, 
\]
from the formula (\ref{daspeed}).
\qed

Since Problem \ref{2is} is an implicit scheme, 
it is not clear whether $\Gamma^{m+1}\in\cO$ can be determined uniquely from the previous polygon 
$\Gamma^m\in\cO$, the time $t_m$, and the time step size $\tau_m$. 
Another question is how to solve the equations
\begin{equation}\label{hjm}
\vech^{m+1}=
\vech^m+\tau_m
F\left(\frac{\Gamma^m+\Gamma^{m+1}}{2}, t_{m+1/2}\right), 
\end{equation}
to obtain (approximation of) $\Gamma^{m+1}$ numerically. 

We fix $\hat{\Gamma}\in\cO$ and $\hat{t}\in [0,T)$, which 
correspond to $\Gamma^m$ and $t_{m+1/2}$, respectively. 
Let $\cK$ be a compact convex set in $\cP^*$ with $\hat{\Gamma}\in\cK\subset\cO$. 
For $\Sigma\in\cK$ and $\hat{\tau}\in(0, \rho(\hat{\Gamma}, \cO)M(\cK, T)^{-1})$, 
we can define $\tilde{\Sigma}\in\cO$ by 
\[
\vech(\tilde{\Sigma})=
\vech(\hat{\Gamma})+\hat{\tau}
F\left(\frac{\Sigma+\hat{\Gamma}}{2}, \hat{t}\right).
\]
In other words, we can define 
$\Lambda(\Sigma):=
\Lambda(\Sigma; \hat{\Gamma}, \hat{t}, \hat{\tau})
:=\tilde{\Sigma}$ which is a map from $\cK$ to $\cO$. 
We have the following lemma.
\begin{Lem}
Let $\eps\in(0, \rho(\hat{\Gamma}, \cO))$ 
and $\lambda\in(0,1)$ be fixed. 
Suppose that $\hat{\tau}$ satisfies the condition
\[
0<\hat{\tau}\leq\min\left\{
T-\hat{t},\ 
\frac{\eps}{M(\hat{\cK}, T)},\ 
\frac{2\lambda}{L(\hat{\cK}, T)}
\right\},
\]
where $\hat{\cK}:=\ov{B(\hat{\Gamma}, \eps)}$. 
Then $\Lambda$ maps $\hat{\cK}$ into $\hat{\cK}$ and 
\begin{equation}\label{dLL}
d(\Lambda(\Sigma^1), \Lambda(\Sigma^2))
\leq \lambda d(\Sigma^1,\Sigma^2)
\quad 
(\Sigma^1, \Sigma^2\in\hat{\cK}). 
\end{equation}
Namely,
$\Lambda$ is a contraction mapping on $\hat{\cK}$ and there exists a unique 
fixed point of $\Lambda$ in $\hat{\cK}$. 
\end{Lem}
\prf 
It is enough to show (\ref{dLL}), which is proved as follows: 
\begin{eqnarray*}
\lefteqn{d(\Lambda(\Sigma^1), \Lambda(\Sigma^2))}\\
&&
=|\vech(\Lambda(\Sigma^1))-\vech(\Lambda(\Sigma^2))|_\infty 
=\hat{\tau}\left|F\left(\frac{\Sigma^1+\hat{\Gamma}}{2},~\hat{t}\right)
-F\left(\frac{\Sigma^2+\hat{\Gamma}}{2}, \hat{t}\right)\right|_\infty \\
&&
\leq \hat{\tau}L(\hat{\cK}, T)
d\left(\frac{\Sigma^1+\hat{\Gamma}}{2}, \frac{\Sigma^2+\hat{\Gamma}}{2}\right) \\
&&
=\hat{\tau}L(\hat{\cK}, T)
\left|\frac{\vech(\Sigma^1)+\vech(\hat{\Gamma})}{2}
-\frac{\vech(\Sigma^2)+\vech(\hat{\Gamma})}{2}\right|_\infty \\
&&
=\frac{\hat{\tau}}{2}L(\hat{\cK}, T)
d\left(\Sigma^1, \Sigma^2\right)
\leq
\lambda d\left(\Sigma^1, \Sigma^2\right).
\end{eqnarray*}
\qed \\
We immediately have the following theorem, 
which gives us efficient numerical scheme to obtain $\Gamma^{m+1}$. 
\begin{Th}\label{th1}
Let $\cK$ be a compact set in $\cO$ and let $\eps\in(0, \rho(\cK, \cO))$. 
We define
\[
\cK_\eps:=\ov{\bigcup_{\Sigma\in\cK}B(\Sigma, \eps)}. 
\]
For fixed $m$ $(<\bar{m})$ in Problem {\rm \ref{2is}}, 
we assume that $\Gamma^m\in\cK$ and 
\[
\tau_m\leq 
\min\left\{
\frac{\eps}{M(\cK_\eps, T)}, 
\frac{2\lambda}{L(\cK_\eps, T)}
\right\},
\]
where $\lambda\in(0,1)$. 
Then there exists uniquely 
$\Gamma^{m+1}\in\ov{B(\Gamma^m, \eps)}$ satisfying {\rm (\ref{hjm})}.

Furthermore, $\Gamma^{m+1}$ is the fixed point of 
the contraction $\Lambda_m:=\Lambda(\cdot\,; \Gamma^m, t_m, \tau_m)$ 
in $\ov{B(\Gamma^m, \eps)}$, and is given by 
the limit of $\Lambda_m^\nu(\Gamma^m)$ as $\nu\to\infty$ 
with the following estimate
\[
d(\Gamma^{m+1}, \Lambda_m^\nu(\Gamma^m))
\leq \lambda^\nu d(\Gamma^{m+1}, \Gamma^m)
\quad 
(\nu\in\N).
\]
\end{Th}
\begin{Th} \label{Th:convergence_implicit_2nd_order_scheme}
We suppose that $\{\Gamma(t)\}_{0\leq t\leq T}$ 
be a $C^{k+1}$-class solution of Problem {\rm \ref{PF}} 
for $k=0, 1$, or $2$. 
There exists $\delta^*>0$, $\tau^*>0$, $C>0$ and 
a non-decreasing function $\omega (a)>0$ with 
\begin{equation}\label{decayomega}
\omega(a)=\left\{
\begin{array}{@{}ll}\disp
o(a^k) & \mbox{if $k=0$ or $1$,}\\[5pt]\disp
O(a^2) & \mbox{if $k=2$,}
\end{array}
\right.
\quad \mbox{as}\ a\downarrow 0,
\end{equation}
such that, if $d(\Gamma^*, \Gamma^0)\leq \delta^*$ and $\tau\leq \tau^*$, 
then $\Gamma^m\in\cO$ $(m=1, 2, \ldots, \bar{m})$ 
are inductively determined and 
\[
\max_{0\leq m\leq \bar{m}}
d(\Gamma(t_m), \Gamma^m)\leq 
\omega(\tau)+C d(\Gamma(0), \Gamma^0), 
\]
holds.
\end{Th}
\prf 
We put 
$\hat{\rho}:=\rho\left(\{\Gamma(t);\ 0\leq t\leq T\}, \cO\right)$, 
and fix $\delta\in(0, \hat{\rho})$ and $\eps\in(0, \hat{\rho}-\delta)$. 
We define
\[
\begin{array}{@{}l}\disp
\cK:=\ov{\bigcup_{0\leq t\leq T}B(\Gamma(t), \delta)}, 
\quad 
\cK_\eps:=\ov{\bigcup_{\Sigma\in\cK}B(\Sigma, \eps)}, \\[5pt]\disp
L:=L(\cK_\eps,T),
\quad 
R(a):=e^{\frac{L}{2}}\left(
1-\frac{aL}{2}
\right)^{-1/a}
\quad 
(0<a<2/L), \\[5pt]\disp
p_m:=|\vece^m|_\infty+\frac{\omega (\tau)}{L}
\quad (m=0, 1, 2, \ldots, \bar{m}), 
\end{array}
\]
where a non-decreasing function $\omega(a)$ $(0<a<T)$, 
which satisfies (\ref{decayomega}), will be defined in (\ref{defomega}) later. 
Since $R(\cdot)$ is an increasing function and $\lim_{a\downarrow 0}R(a)=e^L$, 
there exists $\delta^*>0$ and $\tau^*>0$ such that 
\[
R(\tau^*)\left(
\delta^*+\frac{\omega(\tau^*)}{L}
\right)
\leq \delta,
\quad 
\tau^* <\min\left(
\frac{\eps}{M(\cK_\eps,T)},\,
\frac{2}{L}
\right).
\]

For $m=0, 1, 2, \ldots, \bar{m}-1$, 
we will prove the following inductive conditions: 
\begin{equation}\label{induction}
\Gamma^m\in\cK, \quad p_m\leq R(\tau)^{t_m}p_0
\quad\Rightarrow\quad
{}^\exists\Gamma^{m+1}\in\cK,
\quad p_{m+1}\leq R(\tau)^{t_{m+1}}p_0. 
\end{equation}
The conditions $\Gamma^0\in\cK$ and $p_0\leq R(\tau)^{0}p_0$ for the case $m=0$ 
are obviously satisfied. 

Let us assume the conditions 
$\Gamma^m\in\cK$ and $p_m\leq R(\tau)^{t_m}p_0$ for a fixed $m$. 
Then, from Theorem \ref{th1}, 
there exists $\Gamma^{m+1}$ uniquely in $\ov{B(\Gamma^m, \eps)}\subset\cK_\eps$, 
and we have 
\begin{eqnarray}
&& \hspace{-1truecm}
\vece^{m+1}-\vece^m
=\vech(t_{m+1})-\vech(t_m)-\tau_m V^m
=\tau_m\left\{\vecxi^m+\left(V(t_{m+1/2})-V^m\right)\right\}, 
\label{ej-ej}\\[5pt]
&& \hspace{-1truecm}
\vecxi^m
:=\frac{\vech(t_{m+1})-\vech(t_m)}{\tau_m}-\dot{\vech}(t_{m+1/2}),
\nonumber
\end{eqnarray}
where $V(t):=(V_1(t), V_2(t), \ldots, V_N(t))^{\rm T}$ and 
$V^m:=(V_1^m, V_2^m, \ldots, V_N^m)^{\rm T}$. 

The last term of (\ref{ej-ej}) is estimated as follows. 
Since $\Gamma(t_{m+1/2})\in\cK\subset\cK_\eps$ 
and $\Gamma^{m+1/2}\in\ov{B(\Gamma^m, \eps)}\subset\cK_\eps$, 
we have
\begin{eqnarray}
\lefteqn{\left|V(t_{m+1/2})-V^m\right|_\infty} \nonumber\\[5pt]
&&
=\left|F(\Gamma(t_{m+1/2}), t_{m+1/2})-F(\Gamma^{m+1/2}, t_{m+1/2})\right|_\infty
\leq L d(\Gamma(t_{m+1/2}), \Gamma^{m+1/2}) \nonumber\\[5pt]
&&
=L \left|\vech(t_{m+1/2})-\frac{\vech^m+\vech^{m+1}}{2}\right|_\infty
=L \left|\frac12(\vece^m+\vece^{m+1})-\veczeta^m\right|_\infty, 
\label{Lee}
\end{eqnarray}
where 
\[
\veczeta^m:=\frac{\vech(t_m)+\vech(t_{m+1})}{2}-\vech(t_{m+1/2}). 
\]
Combining (\ref{ej-ej}) and (\ref{Lee}), we obtain
\begin{eqnarray*}
\lefteqn{|\vece^{m+1}|_\infty 
\leq|\vece^m|_\infty+\tau_m \,|\vecxi^m|_\infty
+\tau_mL\left|\frac12(\vece^m+\vece^{m+1})-\veczeta^m\right|_\infty} \\
&&
\leq |\vece^m|_\infty+\frac{\tau_mL}{2}\left(|\vece^{m+1}|_\infty+|\vece^m|_\infty\right)
+\tau_m(|\vecxi^m|_\infty+L|\veczeta^m|_\infty).
\end{eqnarray*}
By the Taylor expansion, we can obtain 
an non-decreasing function $\omega(a)$ ($0<a<T$) 
which satisfies the condition (\ref{decayomega}) and the inequality 
\begin{equation}\label{defomega}
|\vecxi^m|_\infty+L|\veczeta^m|_\infty\leq\omega(\tau). 
\end{equation}
Hence, we have
\[
\left(1-\frac{\tau_m L}{2}\right)|\vece^{m+1}|_\infty
\leq\left(1+\frac{\tau_m L}{2}\right)
|\vece^m|_\infty+\tau_m\omega(\tau), 
\]
and this inequality is equivalent to
\[
\left(1-\frac{\tau_m L}{2}\right)p_{m+1}
\leq\left(1+\frac{\tau_m L}{2}\right)p_m. 
\]
From the inequalities
\[
\left(1-\frac{\tau_m L}{2}\right)\geq
\left(1-\frac{\tau L}{2}\right)^{\tau_m/\tau}
\quad\mbox{and}\quad
\left(1+\frac{\tau_m L}{2}\right)
\leq e^{\tau_m L/2}, 
\]
we obtain
\[
p_{m+1}
\leq
\left(1-\frac{\tau_m L}{2}\right)^{-1}
\left(1+\frac{\tau_m L}{2}\right)p_m
\leq
R(\tau)^{\tau_m}(R(\tau)^{t_m}p_0)
\leq
R(\tau)^{t_{m+1}}p_0. 
\]
The condition $\Gamma^{m+1}\in \cK$ follows from this estimate as 
\[
|\vece^{m+1}|_\infty
\leq
p_{m+1}
\leq
R(\tau^*)^{t_{m+1}}\,p_0
\leq
R(\tau^*)^T
\left(\delta^*+\frac{\omega(\tau^*)}{L}\right)
\leq
\delta.
\]
Hence, we have proved (\ref{induction}), which 
leads us to the estimate: 
\[
|\vece^m|_\infty
\leq
R(\tau^*)^T
\left(|\vece^0|_\infty+\frac{\omega(\tau)}{L}\right)
-\frac{\omega(\tau)}{L}
\leq
R(\tau^*)^T |\vece^0|_\infty
+
\frac{R(\tau^*)^T-1}{L}\omega(\tau).
\]
The assertion of the theorem is obtained by putting $C:=R(\tau^*)^T$ 
and denoting the above term $L^{-1}(R(\tau^*)^T-1)\omega(\tau)$ 
again by $\omega(\tau)$. 
\qed

\bigskip
\noindent
{\bf Acknowledgement.}\ 
The authors are supported by Czech Technical University in Prague, 
Faculty of Nuclear Sciences and Physical Engineering within the 
Jind\v rich Ne\v cas Center for 
Mathematical Modeling (Project of the Czech Ministry of Education, 
Youth and Sports LC 06052).



\noindent
{\bf Authors:}
\vskip 2mm
\noindent
Michal Bene\v s \\
Faculty of Nuclear Sciences and Physical Engineering \\
Czech Technical University in Prague \\
Trojanova 13, 120 00 Prague, Czech Republic \\
e-mail: benes@kmlinux.fjfi.cvut.cz
\vskip 2mm
\noindent
Masato Kimura \\
Czech Technical University in Prague / \\
Faculty of Mathematics, 
Kyushu University \\
6-10-1 Hakozaki, 
Fukuoka 812-8581, Japan \\
e-mail: masato@math.kyushu-u.ac.jp
\vskip 2mm
\noindent
Shigetoshi Yazaki \\
Czech Technical University in Prague / \\
Faculty of Engineering, 
University of Miyazaki \\
1-1 Gakuen Kibanadai Nishi, 
Miyazaki 889-2192, Japan \\ 
e-mail: yazaki@cc.miyazaki-u.ac.jp

\end{document}